\documentclass[12pt]{amsart}
\makeatletter
\@namedef{subjclassname@2020}{\textup{2020} Mathematics Subject Classification}
\makeatother

\usepackage{amssymb}
\usepackage{float}
\usepackage[normalem]{ulem}
\usepackage{comment}

\usepackage{pgfplots}\pgfplotsset{compat=1.9}
\usepgfplotslibrary{fillbetween}

\usepackage{enumerate}
\usepackage{comment}

\usepackage{caption}
\usepackage{subcaption}

\usepackage[colorlinks=true,linkcolor=blue,citecolor=magenta,urlcolor=cyan]{hyperref}

\theoremstyle{plain}

\theoremstyle{definition}

\newcommand{\eps}{\varepsilon}
\renewcommand{\phi}{\varphi}

\usepackage{graphicx}

\begin{document}

\title{Strange bifurcation diagrams}

\begin{abstract}
We investigate a family of one dimensional maps for which the
bifurcation diagram looks differently than the usual ones. We describe
and exemplify various unique and interesting phenomena arising for this family of
maps.
\end{abstract}

\author
[J.~Bielawski]
{Jakub Bielawski$^1$}
\address{$^1$\,%
Department of Mathematics, Krakow University of Economics, Rakowicka
27, 31-510 Krak\'{o}w, Poland}
\email{bielawsj@uek.krakow.pl}

\author
[T.~Chotibut]
{Thiparat Chotibut$^2$}
\address{$^{2}$\,%
Chula Intelligent and Complex Systems Lab, Department of Physics,
Chulalongkorn University, Bangkok 10330, Thailand}
\email{thiparat.c@chula.ac.th}

\author
[F.~Falniowski]
{Fryderyk Falniowski$^{3}$}
\address{$^{3}$\,%
Department of Mathematics, Krakow University of Economics, Rakowicka
27, 31-510 Krak\'{o}w, Poland}
\email{falniowf@uek.krakow.pl}

\author
[M.~Misiurewicz]
{Micha{\l} Misiurewicz$^{4}$}
\address{$^{4}$\,%
Department of Mathematical Sciences, Indiana University Indianapolis,
402 N. Blackford Street, Indianapolis, Indiana 46202, USA}
\email{mmisiure@iu.edu}

\author
[G.~Piliouras]
{Georgios Piliouras$^{5}$}
\address{$^{5}$\,%
Google DeepMind, London EC4A 3TW, United Kingdom}
\email{gpil@google.com}


\maketitle


\section{Introduction}
When investigating the dynamics of one parameter families of interval
maps, one usually looks at their \emph{bifurcation diagram}
(see~\cite{Devaney}). On the horizontal axis there is the parameter,
while on the vertical axis there are the positions of the points of
the trajectory of a chosen starting point. Of course we cannot plot
the whole trajectory, and we are interested in its eventual behavior.
Therefore, for a given parameter, we skip plotting the first long part
of the trajectory, and then plot a much shorter part of it. For maps
with "good" behavior, that is those with negative Schwarzian
derivative\footnote{Schwarzian derivative of a map is given
  by \[Sf(x)=\frac{f'''(x)}{f'(x)}-\frac 32
  \left(\frac{f''(z)}{f'(z)}\right)^2.\]}, to detect all attractors of
the system it is sufficient to look at trajectories of critical points
of the map.\footnote{To complete the information, in the figures showing
bifurcation diagrams presented in this article, horizontal green lines show the position of
critical points of the map.}

\begin{figure}[h]
\begin{center}
\includegraphics[width=150truemm]{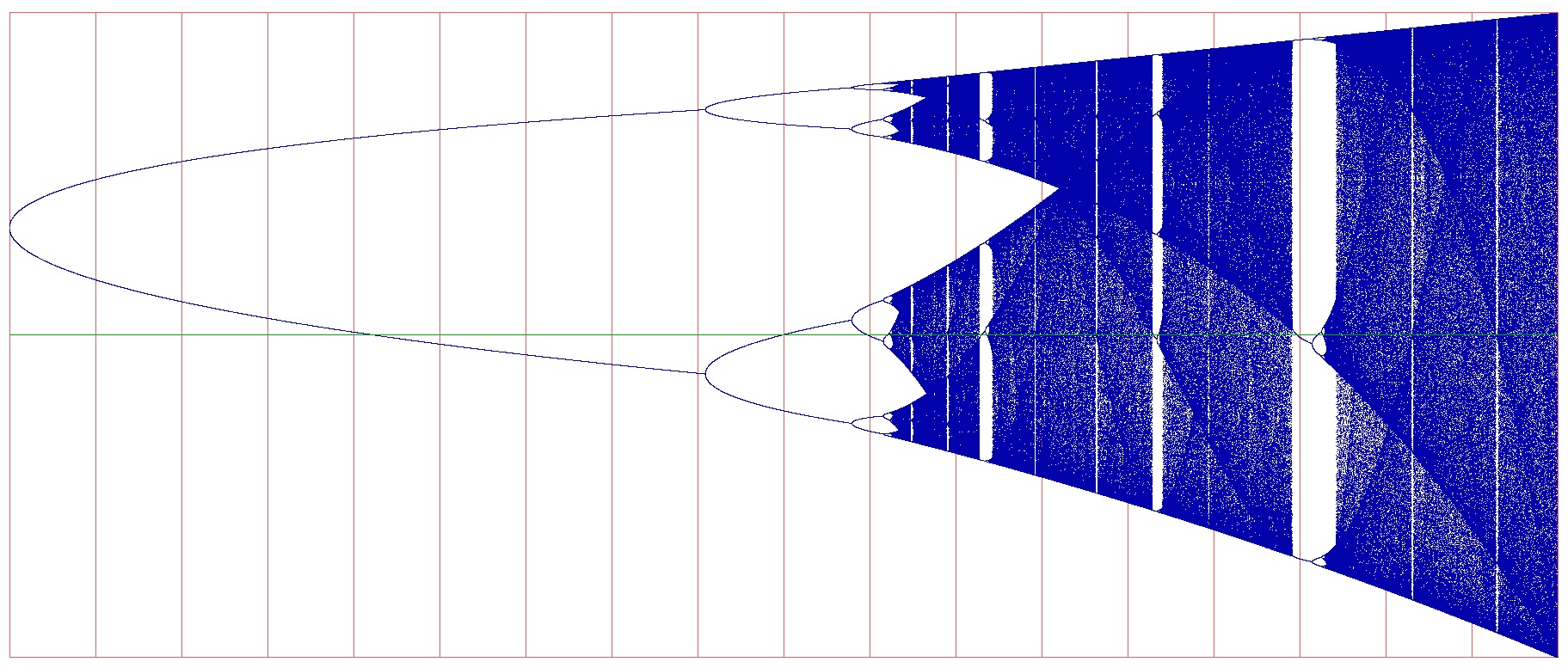}
\caption{Bifurcation diagram for the family of logistic
maps with $r\in [3,4]$.}\label{quad}
\end{center}
\end{figure}

The most popular family is the family of \emph{logistic maps} \cite{Devaney, May}
given by the equation
\begin{equation}
\ell(x)=rx(1-x)
\end{equation}
on the interval $[0,1]$, where $r\in (0,4)$ is the parameter. This is a unimodal map with a critical point at $x_0=1/2$. Its
bifurcation diagram for $r$ from 3 to 4 is shown in Figure~\ref{quad}.

\begin{figure}[h]
\begin{center}
\includegraphics[width=150truemm]{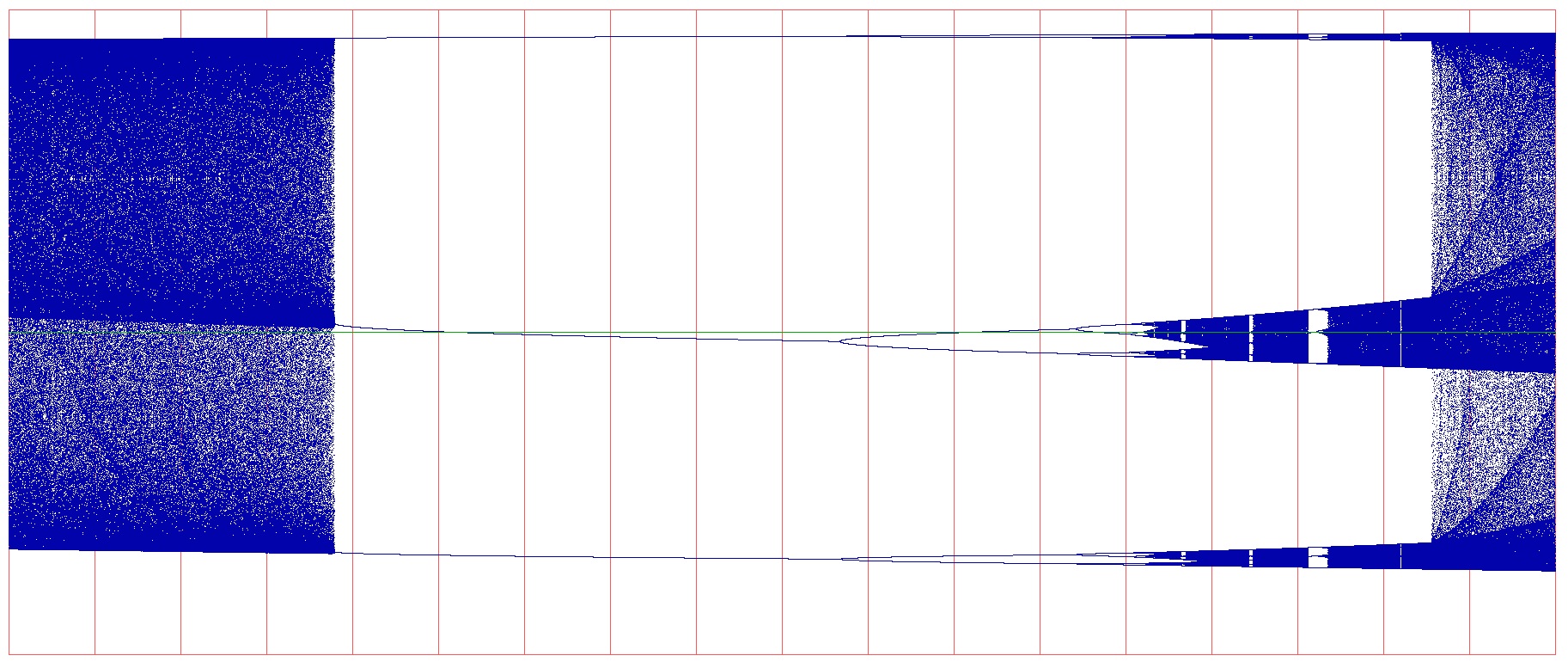}
\caption{Period 3 window for the family of logistic maps.}\label{quad1}
\end{center}
\end{figure}

An important feature of the bifurcation diagram are \emph{periodic
  windows}. They are sets of parameter values for which there exists
an attracting periodic orbit or an attracting cycle of intervals. For
logistic maps, the largest periodic window (excluding initial $2^n$
windows) is the period 3 one, which we see close to $r=3+\frac{5}{6}$.
If we zoom on it, we get a situation illustrated at
Figure~\ref{quad1}.

For other ``natural'' families in most cases the pictures look
similarly. Yet we found a family of maps, appearing in various application contexts,
for which the bifurcation diagram in periodic windows looks
differently. Namely, for some values of parameter, small changes of
the parameter cause large changes of the position of the attracting
periodic points. This family is the family of \emph{EOS maps}
(see~\cite{EOS,bielawski2022memory,eirola1996chaotic}). In this note we will describe and study phenomena arising for these maps.

\section{EOS maps}
The family of \emph{EOS maps}
is defined by
\begin{equation}\label{EOSmap}
F(x)=x+b-\frac1{e^{-ax}+1}
\end{equation}
on the interval $[b-1,b]$, where $a$ is fixed and $b$ is the
parameter. In this note we will look mainly at bifurcation diagrams
for EOS maps with $a=100$. We will call these maps \emph{EOS100 maps}.
Thus, we get the one parameter family of maps
\begin{equation}
F(x)=x+b-\frac1{e^{-100x}+1}.
\end{equation}
We also mention briefly, mainly for descriptive reasons, \emph{EOS200}
maps (the case of $a=200$).

\begin{figure}[h]
\begin{center}
\includegraphics[width=150truemm]{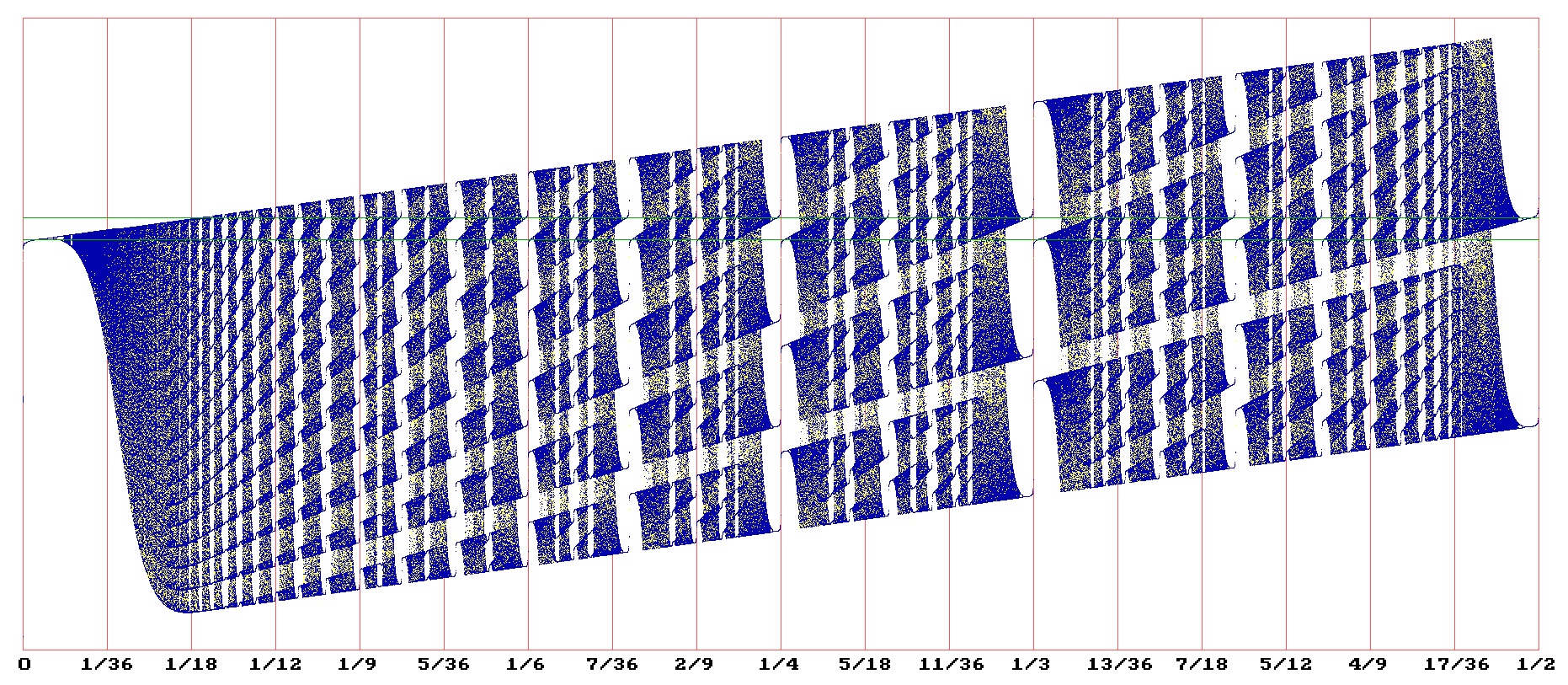}
\caption{Bifurcation diagram for the family of EOS100
  maps.}\label{sbd1a}
\end{center}
\end{figure}

The relevant values of $b$ are $0\le b\le1$, but because of the
symmetry (see~\cite{EOS}), the bifurcation diagram for $1/2\le b\le1$
looks like the bifurcation diagram for $0\le b\le1/2$ turned by 180
degrees. Thus, we show only the part for $0\le b\le1/2$. We draw the
bifurcation diagram for the EOS100 family; see Figure~\ref{sbd1a}. In
the pictures we skip the first 20000 points of the trajectory, and
plot the next 500 ones.

\begin{figure}[h]
\begin{center}
\includegraphics[width=150truemm]{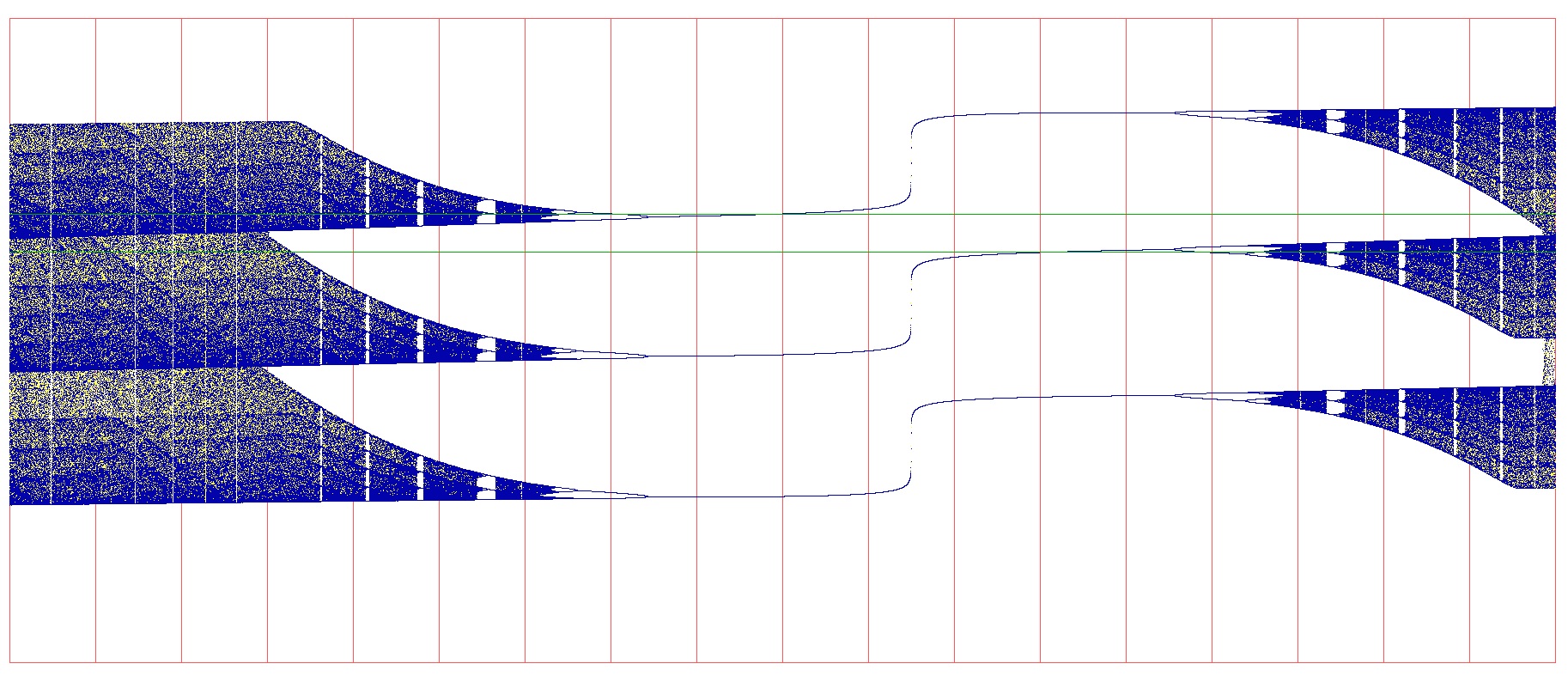}
\caption{Period 3 window for the family of EOS100 maps.}\label{sbd2a}
\end{center}
\end{figure}

The EOS100 maps are bimodal (have two critical points) and have negative Schwar\-zian derivative
\cite{EOS}. Thus, to get full picture it is sufficient to study
trajectories of critical points. We use different colors (yellow and
blue) for the starting point left and right critical point
respectively. We are plotting first yellow points, and then blue, so
many yellow points are hidden below the blue ones. Horizontal green
lines show positions of critical points. We zoom on the period 3
window (Figure~\ref{sbd2a}) to see the finer structure of the
bifurcation diagram.

\section{Hopping attracting periodic orbit}
In Figures~\ref{sbd1a} and \ref{sbd2a} we see that small changes of the parameter result in large changes of the position of the attractor.
The phenomenon that we see can be described as follows. We will
focus on the period 3 window. In this window, choose one of the
period 3 attracting points (smallest, middle or largest). Denote by
$P(b)$ its position. It is easy to see that the function $P$ is
differentiable. Now, while far from $b=1/3$ its derivative is small
(as it should be in the ``normal'' situation), very close to $b=1/3$
this derivative is close to infinity.

\begin{figure}
     \centering
     \begin{subfigure}[b]{0.45\textwidth}
         \centering
         \includegraphics[width=\textwidth]{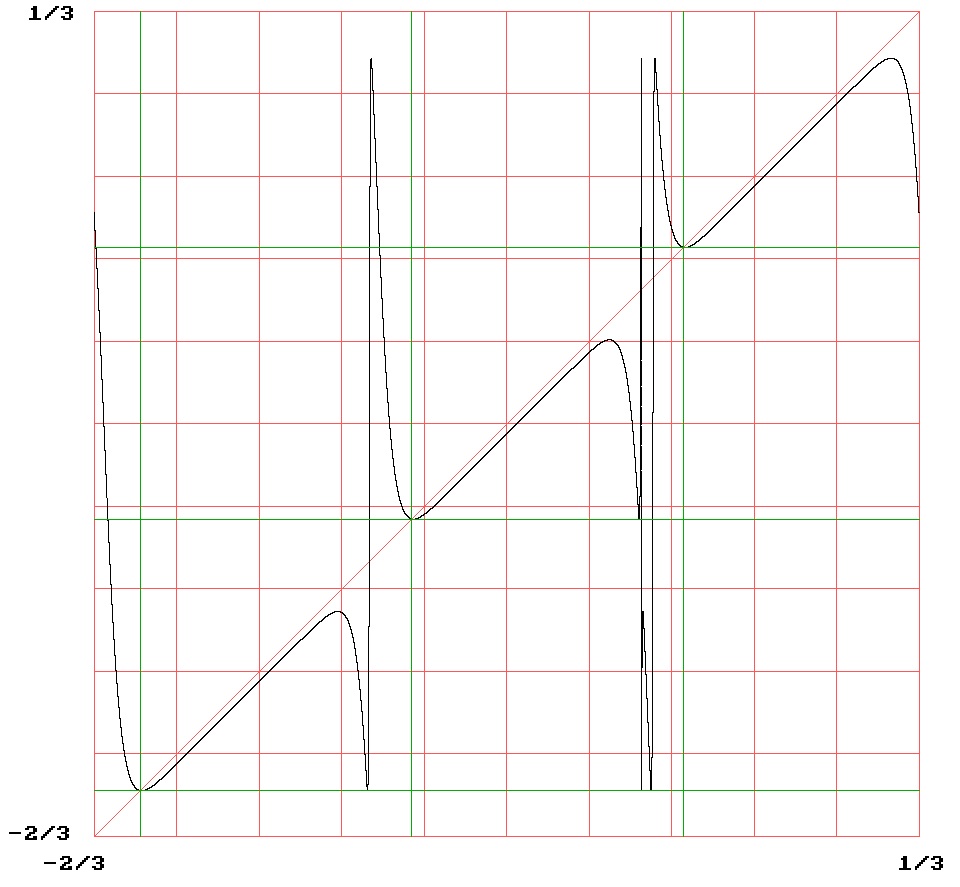}
        \caption{$b=1/3-0.004$}\label{eminus}
     \end{subfigure}
     \hfill
     \begin{subfigure}[b]{0.45\textwidth}
         \centering
         \includegraphics[width=\textwidth]{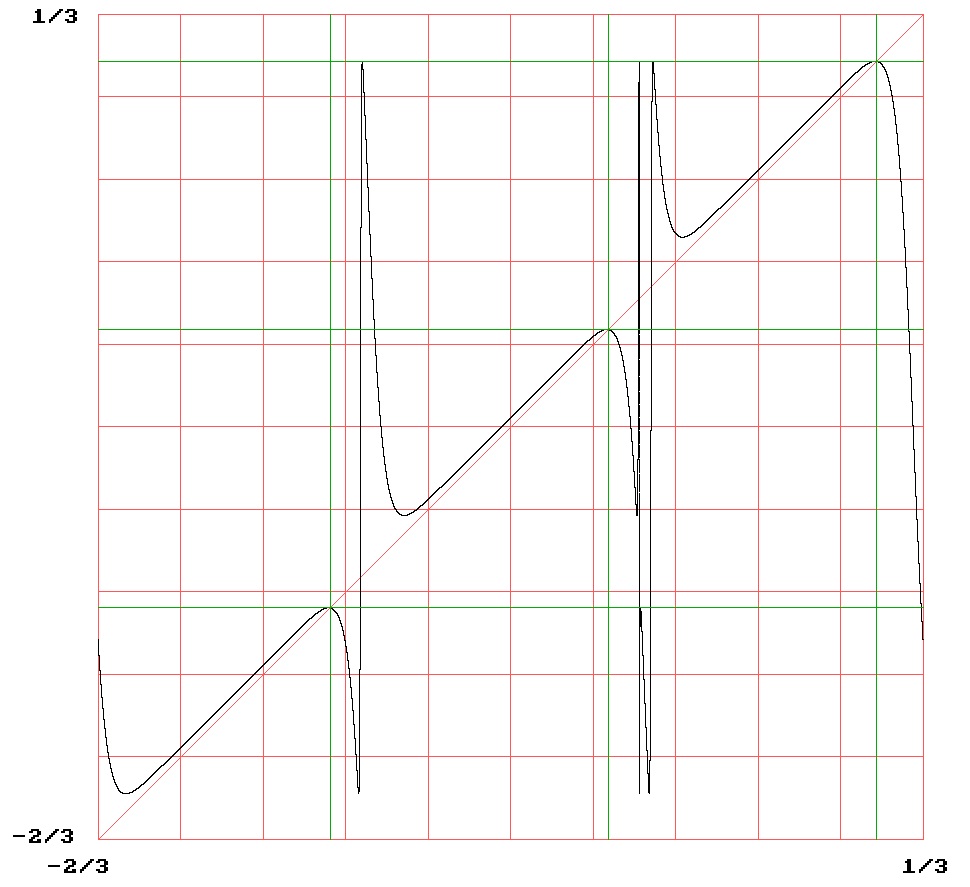}
      \caption{ $b=1/3+0.004$}\label{eplus}
     \end{subfigure}
     \hfill
     \begin{subfigure}[b]{0.45\textwidth}
         \centering
         \includegraphics[width=\textwidth]{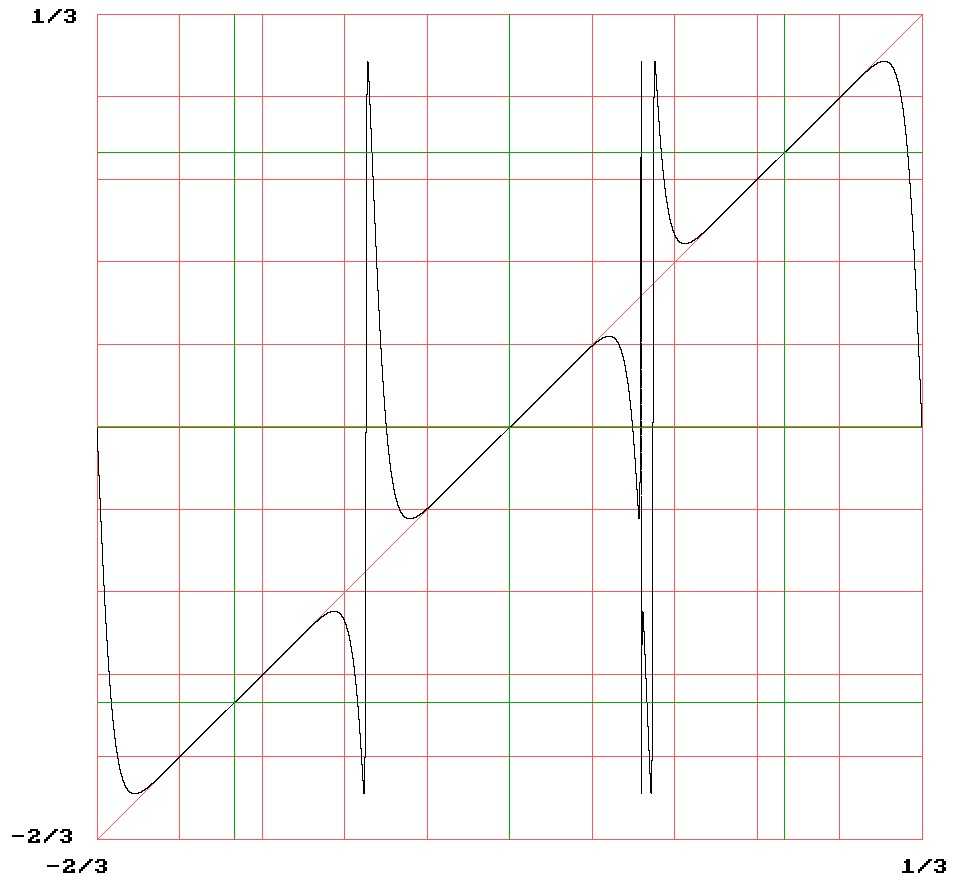}
        \caption{ $b=1/3$}\label{e0}
     \end{subfigure}
        \caption{The graph of $F^3$ for the EOS map for values of $b$ close to $1/3$.}
        \label{fig:three graphs}
\end{figure}

To explain this phenomenon, let us look at the graphs of $F^3$ for
$b=1/3$ and for $b=1/3\pm\eps$ for small $\eps>0$ (see
Figures~\ref{eminus},~\ref{eplus} and ~\ref{e0}). Green horizontal and
vertical lines mark the positions of attracting fixed points of $F^3$.

We see that in Figures~\ref{eminus} and~\ref{eplus} the attracting
fixed points of $F^3$ lie in the region where the graph of $F^3$ looks
more or less like a parabola. This is the ``normal'' behavior, so the
derivative of $P$ is rather small. On the other hand, in
Figure~\ref{e0} the attracting fixed points of $F^3$ lie in the region
where the graph of $F^3$ is almost indistinguishable from the
identity (although it is easy to check that the derivative there is
less than 1). Thus, small changes of $b$ result in the large changes
in $P(b)$.

\begin{figure}[h]
\begin{center}
\includegraphics[width=150truemm]{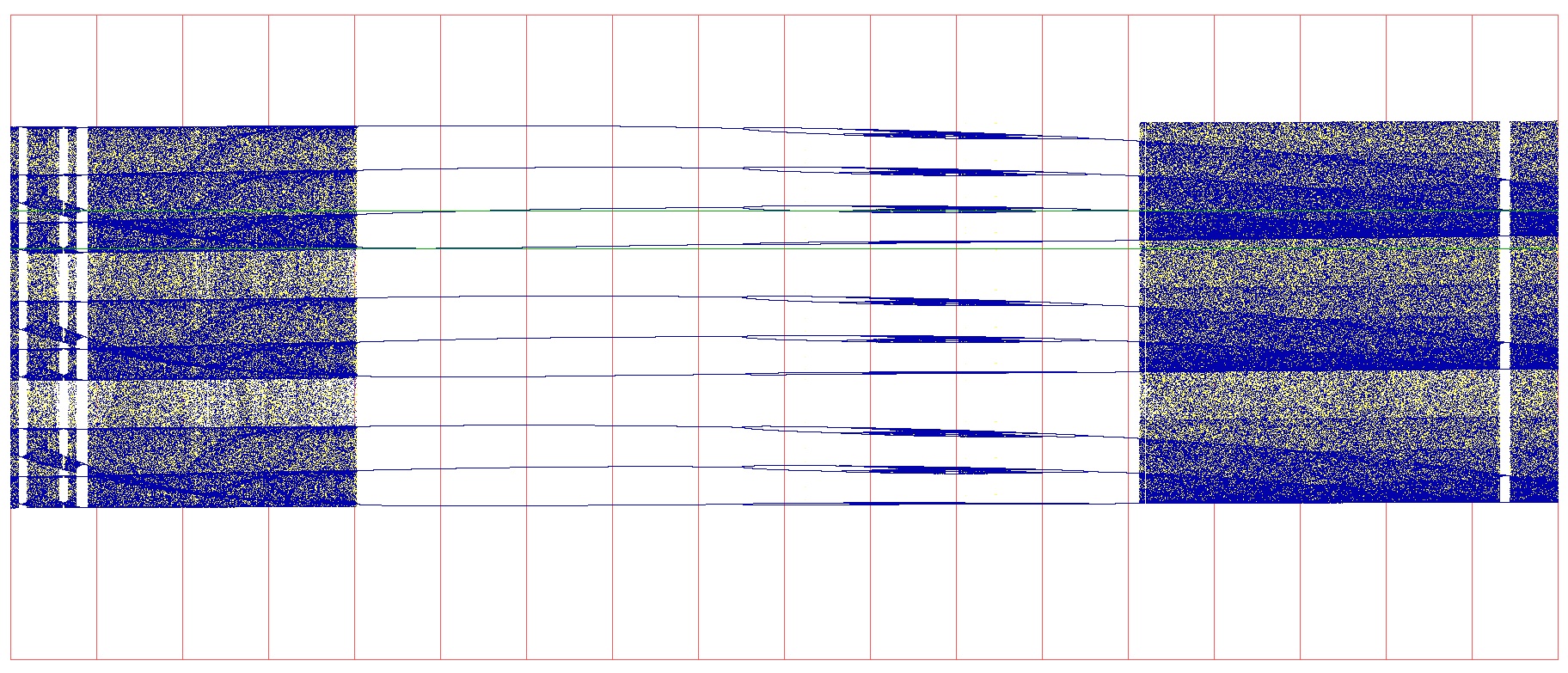}
\caption{Period 10 window for the family of EOS100 maps.}\label{sbd3a}
\end{center}
\end{figure}

\begin{figure}[h]
\begin{center}
\includegraphics[width=100truemm]{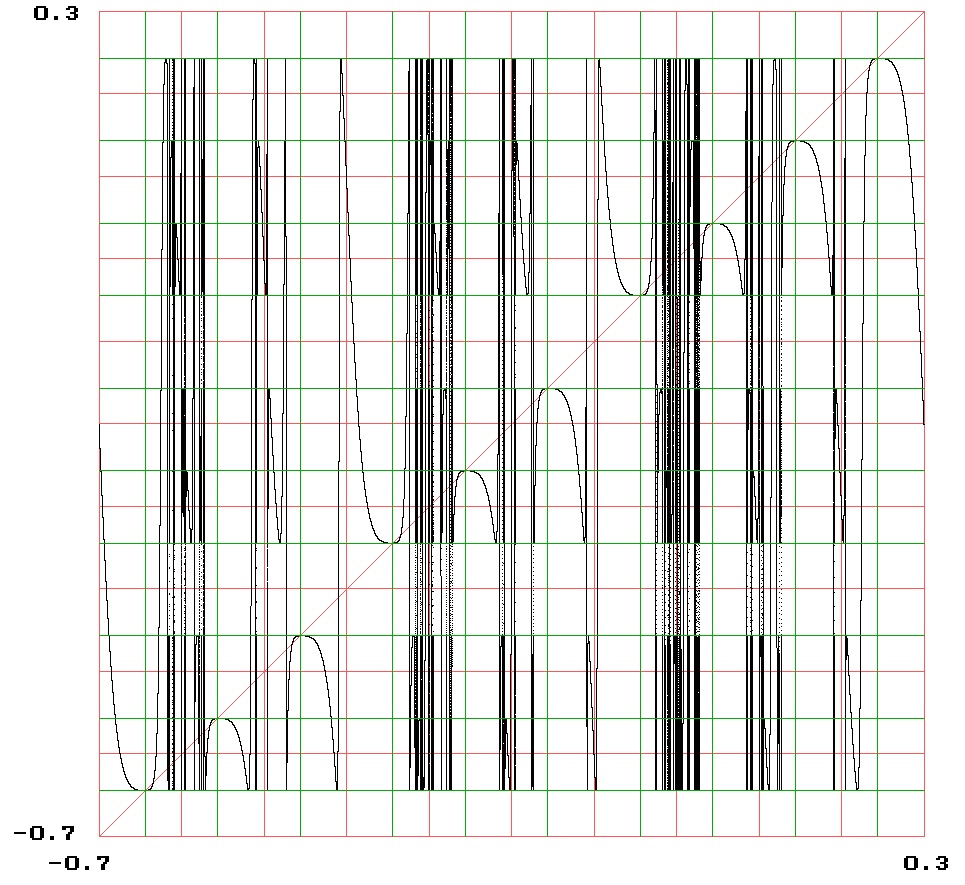}
\caption{The graph of $F^{10}$ for the EOS100 map with
  $b=3/10$.}\label{b03}
\end{center}
\end{figure}

When we look at Figure~\ref{sbd1a}, we may notice that not all
periodic windows look in the way described above. For instance, zooming
at the period 10 window around $b=3/10$ (see Figure~\ref{sbd3a}) , we observe that
 more ``normal'' behavior. Looking at the graph of
$F^{10}$ for $b=3/10$, shown in Figure~\ref{b03}, we see that this iteration does not behave like one described in Figures~\ref{eminus}-\ref{e0}.

\begin{figure}[h]
\begin{center}
\includegraphics[width=150truemm]{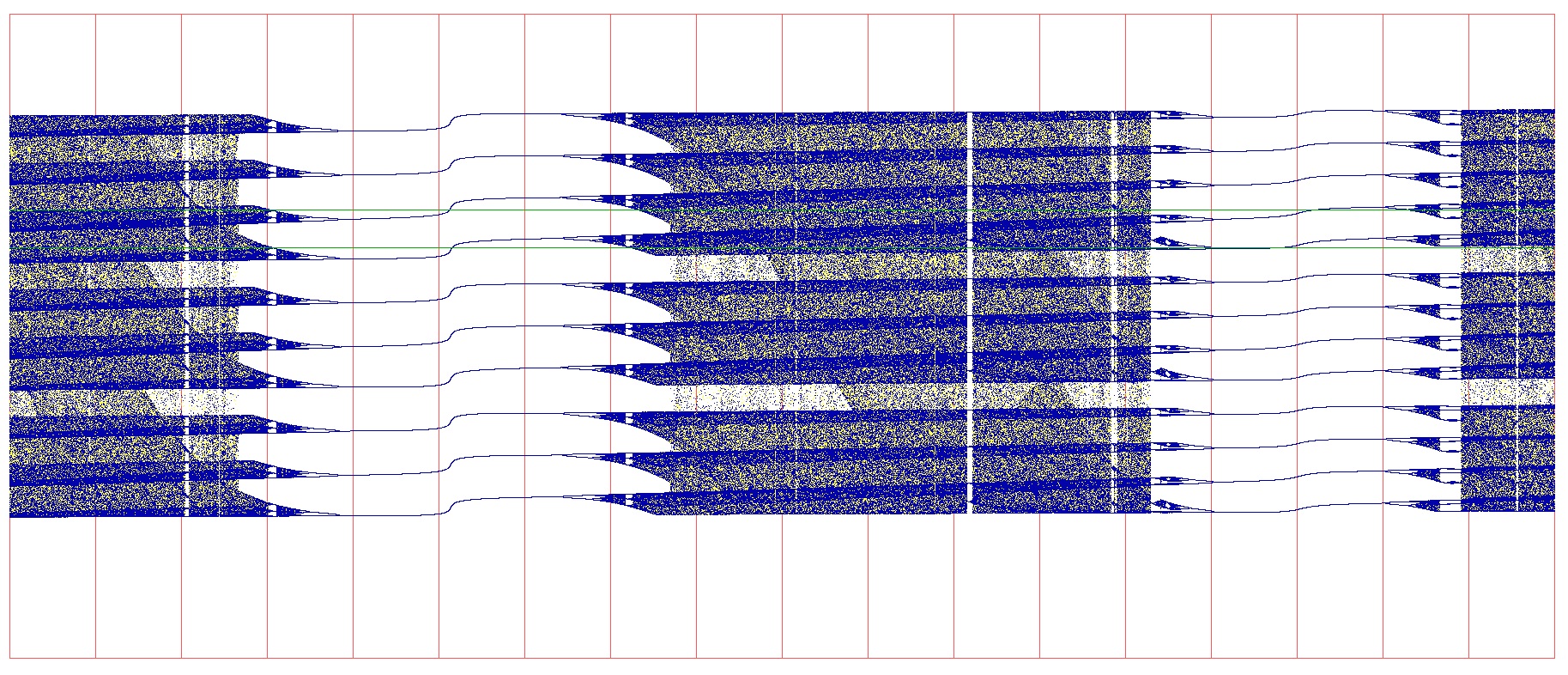}
\caption{Period 10 window for the family of EOS200 maps.}\label{sbd3aa}
\end{center}
\end{figure}

\begin{figure}[h]
\begin{center}
\includegraphics[width=100truemm]{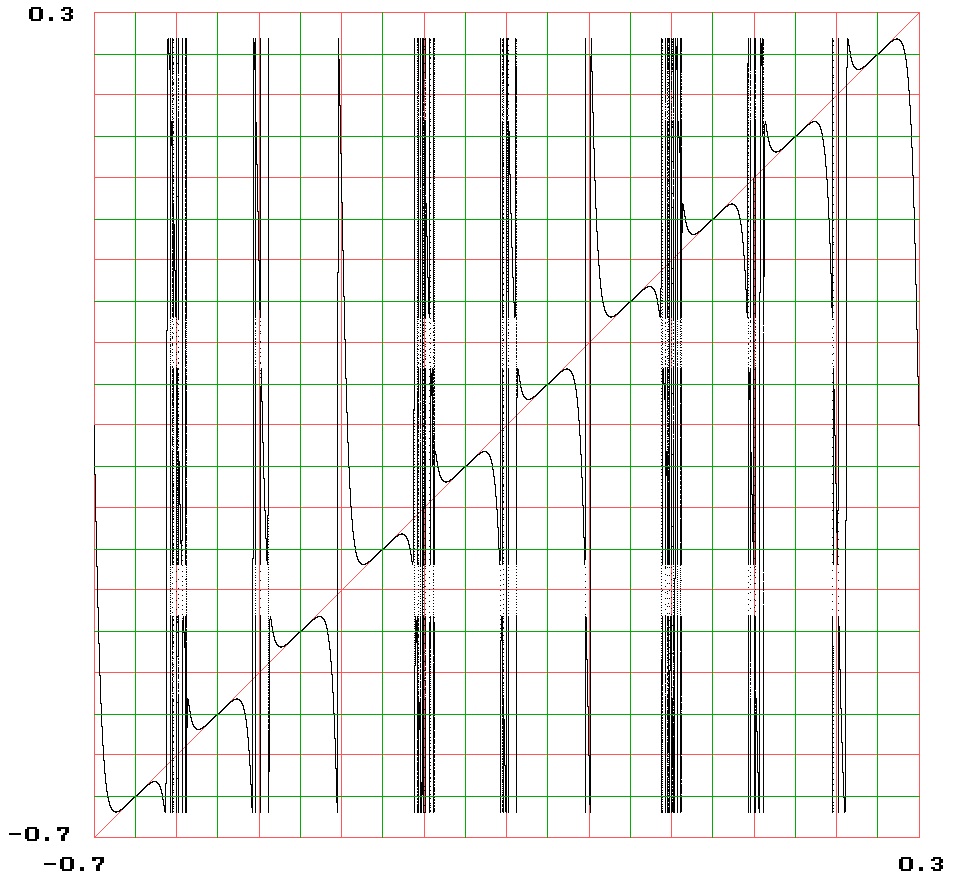}
\caption{The graph of the 10th iterate for $b=3/10$ for the EOS200
  map.}\label{b03b}
\end{center}
\end{figure}

To explain it, let us note that for a window of period $n$ existence
of relatively long intervals where $F^n$ is very close to the identity
occurs for sufficiently large $a$ (see~\cite{EOS}). That may mean that
for this window the value of $a$ is too small. Thus, let us look at
Figures~\ref{sbd3aa} and~\ref{b03b}, made for EOS200 family. There we
see the phenomenon described earlier. To get the derivative of $P$ at
$b=3/10$ even steeper, we should take $a$ even larger. However, if $a$
is too large, then the computational errors become too large, and the
pictures are unreliable. Nevertheless, one should expect appearance of hopping attracting periodic orbit  for each window for sufficiently large values of $a$. 

\section{Different rotational types}
The windows like the one for $b=1/3$ for the EOS100 map are of the
\emph{first rotational type}, that is, for $b$ in some subinterval of
the window, the order of the consecutive points of the trajectory is
the same as for a rotation of the circle, and all points of the orbit
belong to the first or third lap (where the map is increasing).

\begin{figure}[h]
\begin{center}
\includegraphics[width=150truemm]{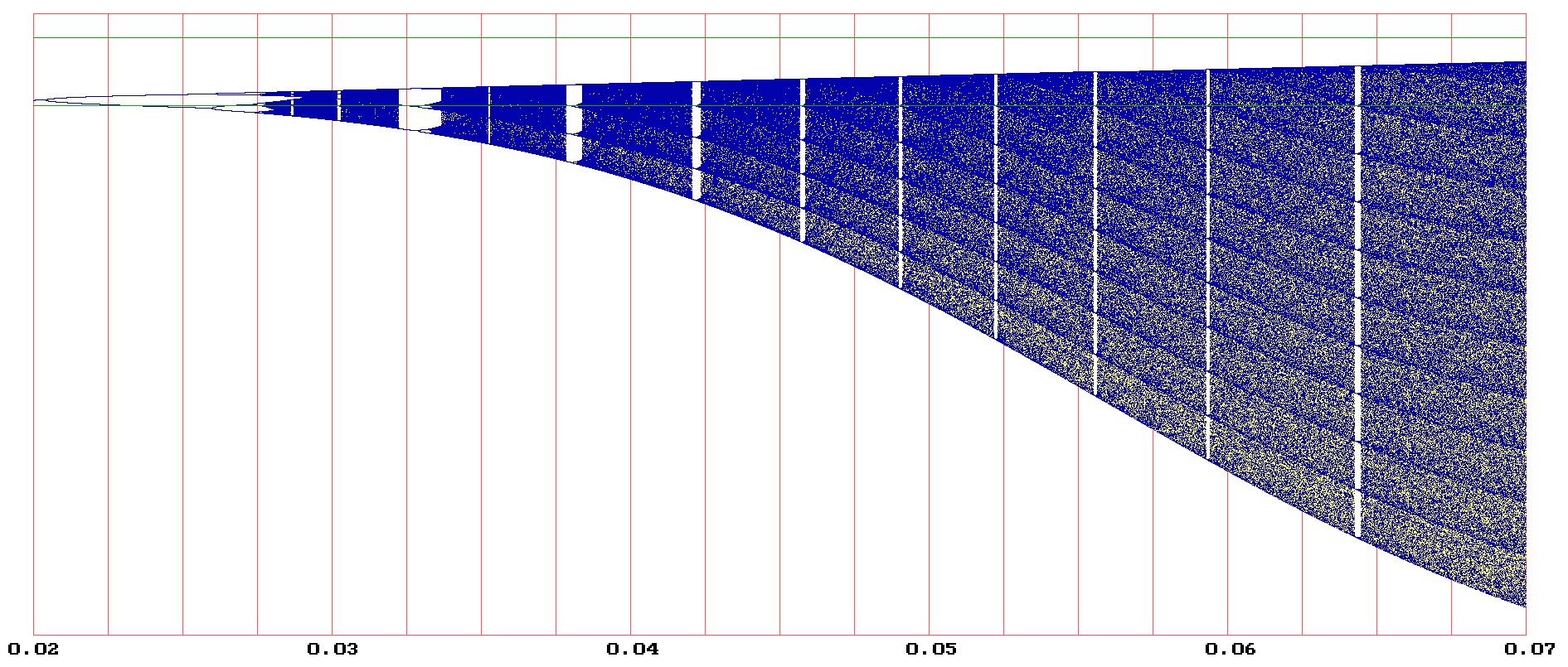}
\caption{Zoom of the bifurcation diagram for the family of EOS100
  maps.}\label{sbd1e}
\end{center}
\end{figure}

Of course there are many more windows, where one or more points of the
attracting periodic orbit are in the second lap. Among those there are
also windows of \emph{second rotational type}, that is, again the
order of the consecutive points of the trajectory is the same as for a
rotation of the circle, but one point of the orbit belongs to the
second lap. 
appear for small values of $b$, in the order of the period of the
attracting periodic point. We can see them as prominent windows in the
zoom of the part of the bifurcation diagram (Figure~\ref{sbd1e}).

\begin{figure}[h]
\begin{center}
\includegraphics[width=150truemm]{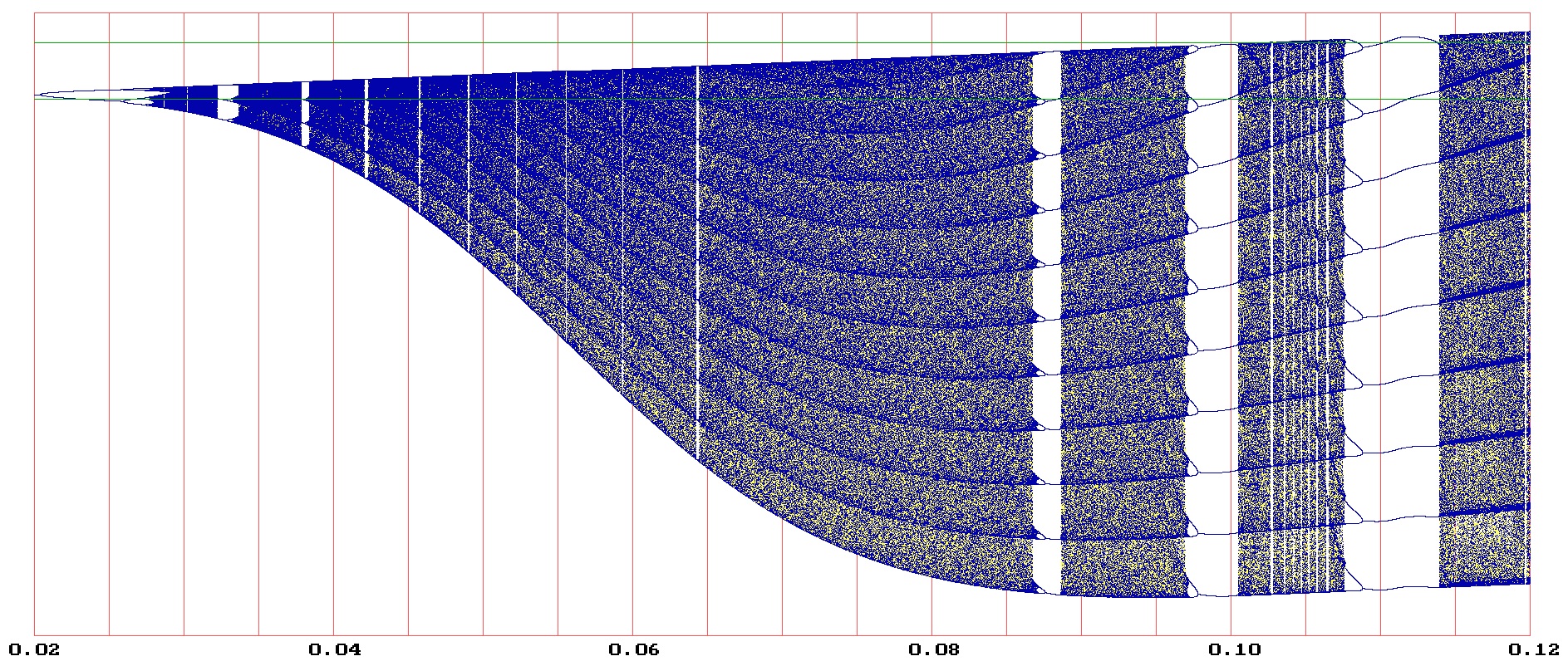}
\caption{Another zoom of the bifurcation diagram for the family of
  EOS100 maps.}\label{sbd1f}
\end{center}
\end{figure}

The order of periodic windows for small values of $b$ is the same as
for the logistic map. The reason is that for small values of $b$ there
is a globally attracting interval, on which the map is unimodal.
However, as $b$ grows, the situation changes, and we see a bimodal map
in the invariant interval. Then the order of periodic windows changes.
Figure~\ref{sbd1f} suggests that for the EOS100 maps the windows of
the second rotational type persist up to period 11, when they are
replaced by the windows of the first rotational type, starting with
the rotation number $1/10$ (at least for relatively small
denominators). However, we see at Figure~\ref{sbd1f} another large
window of period 11 (around $b=0.088$), which is of rotational type,
but seems to be different from types one and two. It is unclear
whether existence of a similar window occurs for all values of $a$.

\section{Phenomena for small values of $b$}
\begin{figure}[h]
\begin{center}
\includegraphics[width=150truemm]{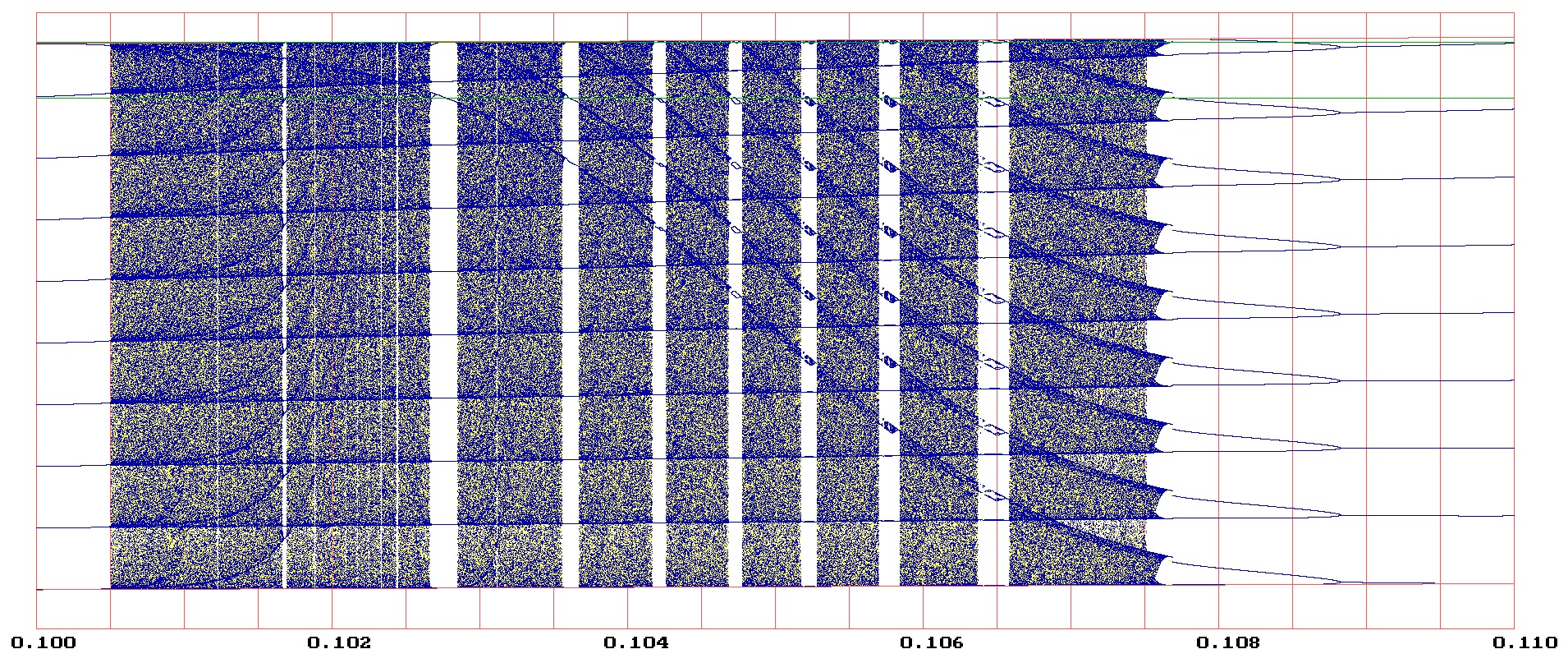}
\caption{Yet another zoom of the bifurcation diagram for the family of
  EOS100 maps.}\label{sbd1i}
\end{center}
\end{figure}

There is another interesting phenomenon for relatively small values of
$b$. There are many windows visible (that is, not too small) between
$b=0.1$ and $b=0.11$. Let us zoom on that region, see
Figure~\ref{sbd1i}.

\begin{figure}[h]
\begin{center}
\includegraphics[width=150truemm]{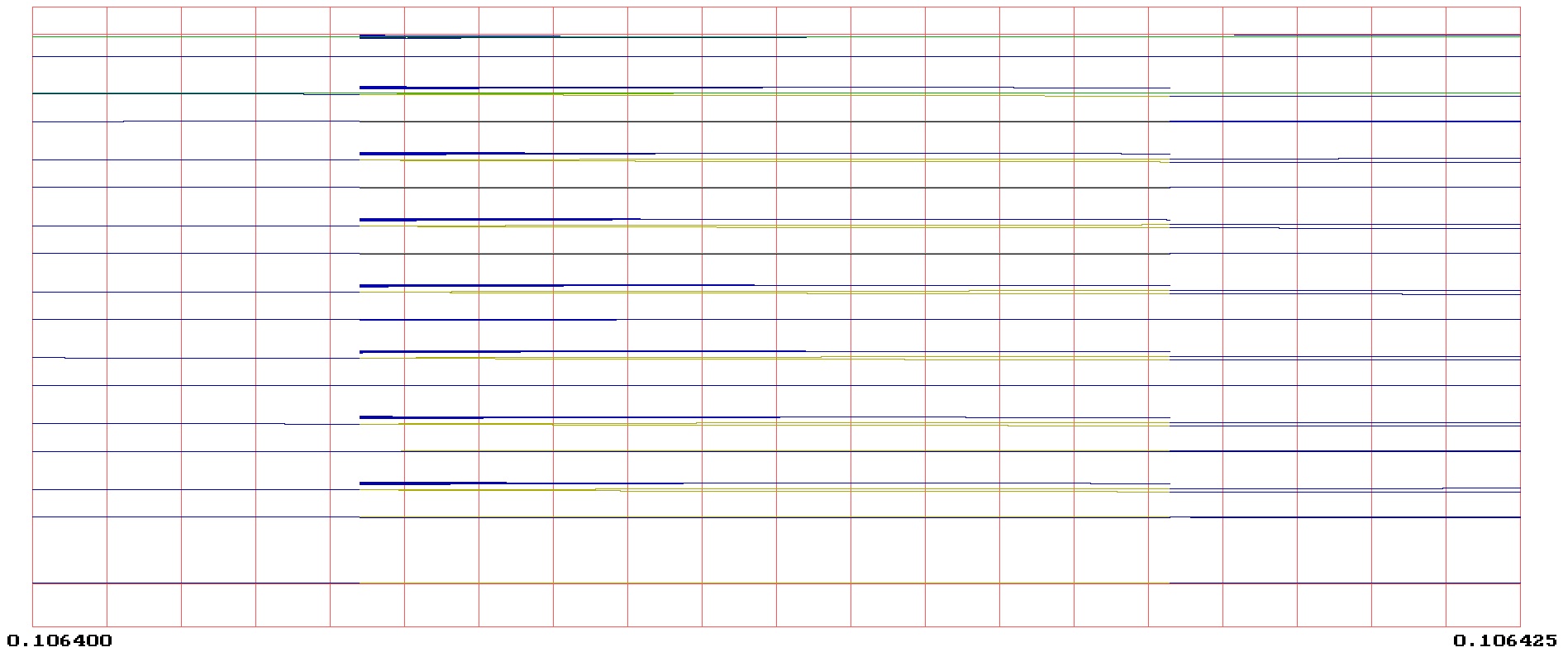}
\caption{Zoom of a part of the period 17 window. Color yellow is
  replaced by olive for better visibility.}\label{sbd1lc}
\end{center}
\end{figure}

We get a quite regular picture, with periods increasing as $b$
increases. The periodic points in those windows are not of rotational
type. Let us now look closer at the period $17$ window. If we look at
a part of this window (see Figure~\ref{sbd1lc}), we see another
interesting phenomenon. Let us recall that the colors yellow and blue
correspond to starting points at two critical points. For the yellow
critical point we see that the periodic point to which its trajectory
converges seems to be a continuous function of $b$. This is not true
for the blue critical point. At about $b=0.1064055$ we see a ``jump''.
Moreover, after this jump the periodic orbit gets replaced by a
periodic cycle of intervals.

This can be explained as follows. The attracting periodic
orbit has the yellow (olive) critical point in its basin of attraction
for the whole time (in the picture, sometimes the yellow line is
covered by the blue one). However, there is also an invariant cycle of
small intervals. The blue critical point is first also in the basin of
attraction of the periodic orbit. However, as $b$ increases, this
critical point enters one of the intervals of the cycle. Thus, we see
a jump.


\begin{thebibliography}{9}


\bibitem{EOS}
J.~Bielawski, T.~Chotibut, F.~Falniowski, M.~Misiurewicz, and
G.~Piliouras.
\newblock {\em
Interval maps mimicking circle rotations},
Communications in Nonlinear Science and
Numerical Simulation {\bf 150}, 108963 (2025)
https://doi.org/10.1016/j.cnsns.2025.108963

\bibitem{bielawski2022memory}
J.~Bielawski, T.~Chotibut, F.~Falniowski, M.~Misiurewicz, and G.~Piliouras.
\newblock Memory loss can prevent chaos in games dynamics.
\newblock {\em Chaos: An Interdisciplinary Journal of Nonlinear Science},
  34(1), 2024.
    {\tt  \url{https://doi.org/10.1063/5.0184318}}

\bibitem{Devaney}
Robert L. Devaney.
\newblock {\em
An Introduction To Chaotic Dynamical Systems, Third Edition
(Addison-Wesley Studies in Nonlinearity)},
Chapman and Hall/CRC 2021.

\bibitem{eirola1996chaotic}
T.~Eirola, A.~V. Osipov, and G.~S{\"o}derbacka.
\newblock {\em Chaotic regimes in a dynamical system of the type many
  predators-one prey}.
\newblock Helsinki University of Technology, Math. Report A368, 1996.
{\tt \url{http://www.math.hut.fi/\%7Eteirola/PS/Pedot.ps}}

\bibitem{May}
Robert M. May. 
\newblock{\em Simple mathematical models with very complicated dynamics}, Nature {\bf 261} (5560): 459-467, 1976.
\end{thebibliography}
\end{document}